\documentclass[11pt]{article} 
\usepackage{amssymb,amsmath} 

\def\R{\mathbb{R}}
\def\Z{\mathbb{Z}}
\def\B{{\cal B}}
\newcommand\vol{\operatorname{vol}} 

\begin{document}
\setlength{\parindent}{0pt}
\setlength{\parskip}{0.4cm}
\bibliographystyle{amsplain} 

\thispagestyle{empty} 

\begin{center}

\Large{\bf The volume of the $10^{\text{th}}$ Birkhoff polytope} 
\normalsize

{\sc Matthias Beck and Dennis Pixton} 

\end{center}

The \emph{$n^{\text{th}}$ Birkhoff polytope} is defined as 
  \[ \B_n = \left\{ \left( \begin{array}{ccc} x_{11} & \cdots & x_{1n} \\ \vdots & & \vdots \\ x_{n1} & \dots & x_{nn} \end{array} \right) \in \R^{n^2} : \ x_{jk} \geq 0 , \begin{array}{l} \sum_j x_{jk} = 1 \text{ for all } 1 \leq k \leq n \\ \sum_k x_{jk} = 1 \text{ for all } 1 \leq j \leq n \end{array} \right\} \ , \] 
often described as the set of all $n \times n$ \emph{doubly stochastic matrices}. 
$\B_n$ is a convex polytope with integer vertices. 
A long-standing open problem is the determination of the relative volume of $\B_n$. 
In \cite{beckpixton} we introduced a method of calculating this volume and used 
it to compute $\vol \B_9$. This note is an update on our progress: with the same 
program, we have now computed 

$\vol \B_{10} = $ 
\scriptsize
\[ \frac{727291284016786420977508457990121862548823260052557333386607889}
{828160860106766855125676318796872729344622463533089422677980721388055739956270293750883504892820848640000000} \ . 
\] 
\normalsize 
%  \begin{multline*}
%\vol \B_{10} =  
%727291284016786420977508457990121862548823260052557333386607889/\\
%8281608601067668551256763187968727293446224635330894226779807213\\
%88055739956270293750883504892820848640000000
%\ . 
%\end{multline*}

We computed this using the spare time on most of the 50 Linux workstations
in the Mathematics department at Binghamton.
The total computation time, scaled to a 1GHz processor, was 6160 days, or
almost 17 years.

As in \cite{beckpixton}, we computed part of the \emph{Ehrhart polynomial} of $\B_{10}$, 
that is, the counting function 
  \[ \# \left( t \B_{10} \cap \Z^{100} \right) \ , \] 
a polynomial in the integer variable $t$. The leading term of this polynomial is 
$\vol \B_{10}/10^9$. For details, as well as the computational tricks which were again used 
in our computation, we refer to our paper \cite{beckpixton} and the accompanying web 
site {\tt www.math.binghamton.edu/dennis/Birkhoff}. 

Some final remarks: 

1. Unless the editors of \emph{Discrete \& Computational Geometry} will allow us to insert 
these new results in the final version of \cite{beckpixton}, this note will solely be 
published on the Mathematics ArXiv ({\tt front.math.ucdavis.edu}). 

2. We will \emph{not} attempt to compute $\vol \B_{11}$ with our current algorithm. 

\bibliographystyle{amsplain}
%\bibliography{bib}
\def\cprime{$'$}
\providecommand{\bysame}{\leavevmode\hbox to3em{\hrulefill}\thinspace}
\providecommand{\MR}{\relax\ifhmode\unskip\space\fi MR }
\providecommand{\MRhref}[2]{%
  \href{http://www.ams.org/mathscinet-getitem?mr=#1}{#2}
}
\providecommand{\href}[2]{#2}

\sc Department of Mathematical Sciences\\
    State University of New York\\
    Binghamton, NY 13902-6000\\
\tt matthias@math.binghamton.edu \\ 
    dennis@math.binghamton.edu \\ 

\end{document}